\documentclass[12pt,reqno,sumlimits]{amsart}   
\usepackage{amssymb,amscd,amsmath,epsfig}
\usepackage{%
  array,
  booktabs,
  dcolumn,
  rotating,
  shortvrb,
  tabularx,
  units,
  url,
}
\newtheorem{thm}{Theorem}[section]
\newtheorem{cor}[thm]{Corollary}
\newtheorem{lem}[thm]{Lemma}
\newtheorem{prop}[thm]{Proposition}
\theoremstyle{definition}
\newtheorem{defn}{Definition}[section]
\newtheorem{rem}{Remark}[section]

\newtheorem{exm}[thm]{Example}

\newcommand{\R}{{\mathbb R}}

\newcommand{\C}{{\mathbb C}}
\newcommand{\D}{{\mathbb D}}
\newcommand{\Z}{{\mathbb Z}}

\newcommand{\N}{{\mathbb N}}

\renewcommand{\to}{\longrightarrow}

\newcommand{\tr}{\operatorname{Tr}}

\newsavebox{\savepar}

\numberwithin{equation}{section}

\newcounter{labelflag} 
\newcommand{\labelon}{\setcounter{labelflag}{1}}
\newcommand{\Label}[1]{
                       \ifnum\thelabelflag=1
                          \ifmmode
                             \makebox[0in][l]{\qquad\fbox{\rm#1}}
                          \else
                             \marginpar{\vspace{0.7\baselineskip}
                                        \hspace{-1.1\textwidth}
                                        \fbox{\rm#1}}
                          \fi
                       \fi
                       \label{#1}
                      }
\labelon

\newcommand{\pdo}{\Psi{\rm DO}}
\newcommand{\Met}{{\rm Met}}
\newcommand{\dvol}{{\rm dvol}}
\newcommand{\cutoffint}{{-\hskip -10pt\int}}
\newcommand{\cutoffinteq}{{-\hskip -12pt\int}}
\newcommand{\e}{{\epsilon}}

\newcommand{\Ci}{{C^\infty}}
\newcommand{\Cl}{{C\ell}}

\newcommand{\fp}{{\rm f.p.}}
\newcommand{\TR}{{\rm TR}}
\newcommand{\res}{{\rm res}}

\begin{document}

\title{ Conformal Anomalies via Canonical
Traces}

\author{S. PAYCHA }

\address{Laboratoire de Math\'ematiques, \\
Universit\'e Blaise Pascal (Clermont II), \\
63177 Aub\`iere Cedex, France\\
E-mail: spaycha@math.univ-bpclermont.fr}

\author{S. ROSENBERG}

\address{Department of Mathematics and Statistics, \\
Boston University, \\
Boston, MA 02215, USA\\
E-mail: sr@math.bu.edu}

\maketitle

\begin{abstract}Using Laurent expansions of canonical traces of holomorphic
  families of
classical pseudodifferential operators, we define functionals on the space of
Riemannian metrics and investigate their conformal properties, thereby giving
a unified description of several conformal invariants and anomalies.
\end{abstract}

\section*{Introduction}

In this paper, we use the Kontsevich-Vishik canonical trace to produce a
series of conformal spectral
invariants (or covariants or anomalies)
associated to conformally covariant
pseudodifferential operators.  Although only one covariant is new, the use of
canonical traces provides a systematic treatment of these covariants.

The search for conformal anomalies is
motivated both by string theory and conformal geometry.  Historically,
the variation of functionals ${\mathcal F}$ on the space of
Riemannian metrics $\Met(M)$ on a closed manifold $M$ under conformal
transformations:
$$g\mapsto e^{2f\,g}, \quad f\in \Ci(M, \R)$$ has been a topic of interest
to both mathematicians and physicists 
going back at least to Hermann Weyl
(see Duff \cite{Du} for a historical review of the physics literature, and
Chang \cite{C} for a survey of recent work in mathematics). 
In physics, the study of conformal invariants underwent a revival in 
the early 1980s with Polyakov's work \cite{Po} 
on the conformal anomaly of bosonic
strings, one of the motivating factors behind the development of
determinant line bundles in mathematics.

The conformal anomaly of a Fr\' echet differentiable map ${\mathcal
F}:\Met(M)\to \C$ at a given (background) metric $g$ is the differential at $0$
of $ {\mathcal F}_g: \Ci(M, \R)\to \C$, ${\mathcal F}_g(f) :={\mathcal
F}(e^{2f}\,g).$ Thus the conformal anomaly in the direction $f$ is
$$\delta_f{\mathcal F}_g:= d{\mathcal F}_g(0).f= 
\frac{d}{dt}\left|_{_{_{_{_{t=0}}}}}\right.
{\mathcal F}(e^{2tf} g).$$
A functional ${\mathcal F}$ is  conformally invariant if $\delta_f{\mathcal
  F}_g=0$ for any Riemannian metric $g$ and any smooth function $f$. 
If $\delta_f {\mathcal F}_g= \int_M f(x) \overline{\delta_x {\mathcal
F}_g}(x) \dvol_g(x)$, then  $\overline{\delta_x {\mathcal F}_g}(x)$ is called the
the local (conformal)
anomaly of ${\mathcal F}_g$ (or equivalently of ${\mathcal F}$ in the
background metric $g$).
A
functional ${\mathcal F}(g, x)$ on $\Met(M)\times M$ is
conformally covariant if, roughly speaking, $\delta_f{\mathcal F}$ does not
depend on derivatives of
$f$ and $g$.

Conformal anomalies arise naturally in quantum field theory.
A conformally invariant classical
action ${\mathcal A}(g)$ in a background metric $g$, for example the string theory
or nonlinear sigma model action, does not usually lead to a
conformally invariant 
effective action ${\mathcal W}(g)$, since the
quantization procedure breaks the conformal invariance and hence 
gives rise to a conformal anomaly. In particular, in string theory 
the conformal invariance persists after quantization only in 
specific critical dimensions. 

From a path integral point of view, the conformal anomaly
of the quantized action is often said to arise from a lack of conformal
invariance of the formal measure on the configuration space of the QFT.
Whatever this means, we can detect the source of the conformal anomaly in the
quantization procedure.  
In order to formally reduce the path
integral to a Gaussian integral, one writes the classical action as a
quadratic expression ${\mathcal A}(g)(\phi)= \langle A_g\phi, \phi\rangle_g$ where
$\phi$ is a field, typically a tensor on $M$, $A_g$ a differential operator on
tensors and $\langle\cdot, \cdot\rangle_g$ the inner product 
induced by  $g$. Because this inner product is not conformally
invariant, 
the conformal invariance of ${\mathcal A}(g)$
usually translates to a conformal covariance of the operator $A_g$.  An
operator $A_g$ is {\it conformally covariant } of bidegree $(a,b)\in \R^2$
if
\begin{equation}\label{eq:introconfop}
 A_{\bar g}= e^{-b\, f} A_g e^{a\,f},
\end{equation}
for $\bar g = e^{2f}g.$  
Thus this first step, which turns a conformally invariant quantity (the
classical action) to a conformally covariant operator,  already
breaks the conformal
invariance.

 The second step in the computation of the path integral uses an Ansatz to
  give a meaning to the formal determinants that arise from the Gaussian
  integration. Mimicing finite dimensional computations, the 
effective action derived from a formal integration over the
  configuration space ${\mathcal C}$ is
 $$e^{-\frac{1}{2}{\mathcal W}(g)}:= \int_{\mathcal C}
 e^{-\frac{1}{2}{\mathcal A}(g)(\phi)}\, {\mathcal D}\phi= {\rm
  ``det{\hbox {''}}}(A_g)^{-\frac{1}{2} }.$$ 
If there were
  a well defined determinant ``det'' on differential operators 
with the usual properties,
  (\ref{eq:introconfop}) would yield
\begin{eqnarray*}
 {\rm ``det{\hbox {''}}}(  A_{e^{2f} g})&=& {\rm ``det{\hbox {''}}}
( e^{-b\, f} A_g e^{a\,f})\\
&=& {\rm ``det{\hbox {''}}}( e^{-b \,f}) \, {\rm ``det{\hbox {''}}}( A_g )\, 
{\rm ``det{\hbox {''}}}(e^{a\,f})\\
&=& {\rm ``det{\hbox {''}}}( e^{(a-b)\,  f}) \, 
{\rm ``det{\hbox {''}}}( A_g ),
\end{eqnarray*}
where $e^{c\, f}$ is treated as a multiplication operator for $c\in \R.$
Hence, even if a ``good'' determinant exists, the effective action ${\mathcal
W}(g)$ would still suffer a conformal anomaly, since $A_g$ is only conformally
covariant:
$$\delta_f{\mathcal W}( g)= 
\delta_f\log {\rm ``det{\hbox {''}}}( A_g)=
\delta_f\log {\rm ``det{\hbox {''}}}( e^{(a-b)\, f})=
(a-b)\,{\rm ``tr{\hbox {''}}}(f),$$ 
where ``tr'' is a hypothetical trace associated to
``det''.  

The $\zeta$-determinant $ {\rm det}_\zeta$ on operators
is used by both physicists and
mathematicans as an Ersatz for the
usual determinant on matrices.
Since the work of Wodzicki and Kontsevich--Vishik, we know the
$\zeta$-determinant has a
multiplicative anomaly, which fortunately does not affect our rather
specific situation. Indeed, the above heuristic derivation holds
(Branson-Orsted \cite{BO}, Parker-Rosenberg \cite{PR}, Rosenberg \cite{R}): 
$$ \delta_f \,\log{\rm det}_\zeta ( A_{ g})=(a-b)\,{\rm
tr}^{A_g}(f),
$$ 
if one replaces 
``tr''$(f)$ with ${\rm tr}^{A_g}(f),$
the finite part in the heat-operator expansion ${\rm tr}(f\,
e^{-\e A_g})$ when $\e\to 0$.  
(Here and whenever the heat operator $e^{-\e A_g}$ appears,
we assume that $A_g$ is elliptic with non-negative leading symbol.)
In summary, 
the regularization
procedures involved in the $\zeta$-determinant and the finite
part of the heat-operator expansion are {\it not} 
responsible for the conformal
anomaly of the effective action ${\mathcal W}(g)$; the conformal anomaly 
appears as soon as one uses
the {\it conformally covariant} operator $A_g$ associated to the
originally {\it conformally invariant} action ${\mathcal A}(g)$. 

These QFT arguments lead to the search for conformally covariant
operators and associated spectral conformal covariants.
There are four types of conformal covariants in the literature, in order of
computational difficulty: 
(i) local
covariants, those that depend only on the metric at a fixed point 
(ii) global invariants which are the integrals of (noncovariant) local
quantities, 
(iii)
global invariants which are not integrals of local expressions, but 
 whose variation in any metric direction is local; (iv) global invariants
 which are not integrals of local expression, and whose variation in conformal
 directions is nonlocal.
All four types have examples associated to spectral $\zeta$- and 
$\eta$-functions, as we
now explain.  

For (i), the residue at $z=1$ of the local zeta function
 $\zeta_{A_g}(z,x)$, which
turns out to be  proportional to the local Wodzicki residue
$\res_x(A_g^{-1})$, is a pointwise conformal
 covariant for a conformally covariant
 operator $A_g$, under certain ellipticity and positivity conditions
 on the operator \cite{PR}.  
(A classical example of a pointwise invariant 
is the length of the Weyl tensor \cite{We}.)
For (ii), the value at $z=0$ of the global $\zeta$-function
 $\zeta_{A_g}(z)$ of a conformally covariant operator $A_g$ is conformally
 invariant, again for certain operators, which may be 
pseudodifferential \cite{PR,R}:
$$ \delta_f \zeta_{ A_{g}}(0)=0. $$ 
In hindsight, this can be predicted by thinking of
$\zeta_{ A_{g}}(0)$  as an Ersatz for 
``tr''$({\rm Id})$ in the heuristic notation above.  It is well known that 
$\zeta_{ A_{g}}(0)$  is the integral of the finite part of the pointwise 
heat kernel of $A_g$
(up to the nonlocal conformally invariant term ${\rm dim\ Ker}(A_g)$).  
When $A_g$ is a differential operator, $\zeta_{A_g}(0) = -\frac{1}{{\rm
    ord}(A_g)} \res(\log A_g)$, so 
the conformal invariance of
$\zeta_{A_g}(0)$ is equivalent to
the conformal invariance 
of the
exotic determinant introduced by Wodzicki for zero
order classical pseudodiffferential operators
and
extended by Scott \cite{Sc} 
to the residue determinant $\det_{\res}(A_g)= e^{\res (\log
A_g)}$ on operators of any order.
This gives another description of $\zeta_{A_g}(0)$ as the integral
   of a local quantity, namely the local Wodzicki residue of the logarithm of
   $A_g.$

Jumping to (iv),
conformal anomalies arising from $\zeta$-determinants of  conformally 
covariant operators vanish in certain cases, for one has \cite{PR,R}
$$
\delta_f \,\zeta^\prime_{  A_{ g}}(0)
= -\delta_f \log {\rm det}_\zeta (  A_{ g}) 
 =(a-b)  \int_M f(x)  a_{n}(A_g,x)\dvol_g(x),$$ 
where as $\e\to 0$
$${\rm tr}\left( e^{-\e A_g}\right)\sim \sum_{j=0}^\infty \left(\int_M 
a_j(A_g,x)\dvol_g(x)\right) \e^{\frac{j-n}{2\, \alpha}},$$ 
for $\alpha = {\rm ord}(A_g)$, $n = {\rm dim}(M);$ here we assume
$A_g$ has all but
finitely many eigenvalues nonnegative.
The nonlocal nature of the functional determinant and its variation
is well known; however, the above formula shows it gives rise to a
local conformal anomaly $(a-b) a_n(A_g, x) $.
In particular, 
$ \zeta_{A_g}^\prime(0)$
yields a conformal invariant in odd dimensions, 
as $ a_{n}(A_g)$ then vanishes.
The conformal anomaly $
\delta_f \log {\rm det}_\zeta ( \Delta_{ g})$, where $\Delta_g$ is the 
Laplace-Beltrami operator on a closed Riemannian surface,
is responsible for the
conformal anomaly in bosonic string theory; since the coefficients $a, b$
depend on the dimension of the manifold and the rank of auxiliary tensor
bundles, 
combinations of such conformal
anomalies cancel in certain critical dimensions, viz.~the cancellation 
of conformal anomalies in 26 dimensions for bosonic
string theory \cite{Po}.  Further work on the conformal anomaly of
functional determinants is in work of Branson and Orsted \cite{B1,B2,BO}.

For (iii), if $A_g$ is a self-adjoint invertible elliptic operator, the phase
of its $\zeta$-determinant can be expressed in terms of the $\eta$-invariant
$\eta_{A_g}(0)$ by
$${\rm det}_\zeta(A_g):= {\rm det}_\zeta(\vert A_g\vert)\cdot
e^{i\frac{\pi}{2}\left( \zeta_{\vert A_g\vert} (0)- \eta_{A_g}(0)\right)}.$$
Again, only in certain dimensions is the phase
conformally invariant; namely if ${\rm dim}(M)$ and ${\rm ord}(A_g)$ have
opposite parity \cite{R}.

We will study these four types of
conformal anomalies and covariants in the
common framework of variations of Kontsevich-Vishik
functionals of conformally covariant operators. 
Whereas previous work on
conformal anomalies uses heat kernel expansions, we use
$\zeta$-function techniques instead. Our starting point is canonical traces,
which are cut-off
integrals of symbols of non-integer order pseudodifferential operators,  
which extend
 to Laurent
expansions of cut-off integrals of holomorphic families of
symbols. 
These coefficients are universal expressions in
the symbol expansion of the family (Paycha-Scott \cite{PS}), 
so their regularity
properties and their variation in terms of external parameters (here the
metric) are easily controlled.  We thereby avoid some 
technical difficulties in the
variation of
heat kernel asymptotic expansions. The main result of the paper is that
the coefficients of the Laurent expansions
give explicit conformal
anomalies. 

In more detail,
the three functionals
$\zeta_{A_g}(0)$, $\zeta^\prime_{A_g}(0)$ and $\eta_{A_g}(0)$ are all
$A_g$-weighted traces in the notation of \cite{P2},
namely ${\rm tr}^{A_g}(I)$, ${\rm tr}^{A_g}(\log A_g)$
and ${\rm tr}^{ A_g }(A_g\, \vert A_g\vert^{-1})$ respectively.  
Here, for a weight $Q$ (i.e.~an admissible positive order elliptic operator), 
the $Q$-weighted trace $\tr^Q(A)$ of a classical
pseudodifferential operator $A$ is the finite part at $z=0$ of the
meromorphic map $z\mapsto {\rm TR}(A \, Q^{-z})$ (up to a factor depending on
the kernel of $Q$), where ${\rm TR}$ is the Kontsevich-Vishik
canonical trace on noninteger order operators extending the usual trace on
smoothing operators \cite{KV}. (This definition of weighted trace is
equivalent to previous ones \cite{P2} by the discussion after Def.~3.)
Thus all our spectral invariants are examples of canonical traces.

If the conformally covariant
operator $A_g$ is a weight, we may define functionals given by meromorphic
functions $z\mapsto {\mathcal F}_h(g)(z):= {\rm TR}(h(A_g) \, A_g^{-z})$ where $h$
is a real or complex valued function defined on a subset $W\subset \C$. In
particular, the
functionals $\zeta_{A_g}(z)$ and $\eta_{A_g}(z)$ 
 correspond to choosing $h(\lambda)=1$ (with $W=\C$) and
$h(\lambda)=\frac{ \lambda}{\vert \lambda\vert}$ (with $W=\R/\{0\}$). Using
results
on the coefficients in the Laurent expansion \cite{PS} for
$z\mapsto{\mathcal F}_h(g)(z)$ at $z=0$,
we derive the conformal anomaly of these meromorphic
functionals (Theorem \ref{thm:confanom}):
$$
\delta_f {\rm TR}( h( A_{ g})\, A_{ g}^{-z})
= (a-b)  {\rm TR}\left(f\, h^\prime( A_{ g})\,A_{g}^{-z+1}\right)
-z\,(a-b) {\rm TR}\left(f\, h( A_{ g})\,A_{g}^{-z}\right).
$$
This formula strongly depends on the tracial nature
 of the canonical trace TR on
noninteger order operators (\ref{eq:TRcyclic}).
Identifying the coefficients on either side, we get a hierarchy of functionals
and their conformal anomalies, the first one involving the Wodzicki residue 
${\rm res}$:
\begin{eqnarray*}
 \delta_f {\rm res}( h( A_{ g})) 
& =&(a-b)\,  {\rm res}\left(f\, h^\prime( A_{ g})\,A_{g}\right);\\
 \delta_f {\rm tr}^{A_{ g}}( h( A_{ g})) 
 &=&(a-b) \, {\rm tr}^{A_g}\left(f\, h^\prime( A_{ g})\,A_{g}\right)
+\frac{a-b}{\alpha}\,{\rm res}\left(f\, h( A_{ g})\right);\\
 \delta_f {\rm tr}^{A_{g}}( h( A_{g})\,\log A_g)) 
 &=&(a-b) \, {\rm tr}^{A_g}\left(f\, h^\prime( A_{ g})\,A_{g}\, \log
 A_g\right)\\
&&\qquad +\frac{b-a}{\alpha}\,{\rm tr}^{A_g}\left(f\, h( A_{ g})\right);\\ 
&\vdots&\\
 \delta_f{\rm tr}^{A_{ g}}( h( A_{g})\,\log^j A_g)) 
 &=&(a-b)\,  {\rm tr}^{A_g}\left(f\, h^\prime( A_{ g})\,A_{g}\, \log^j
 A_g\right)\\
&&\qquad +j\,\frac{a-b}{\alpha}\,{\rm tr}^{A_g}
\left(f\, h( A_{ g})\, \log^{j-1} A_g\right).\\
\end{eqnarray*}
Different choices for $h$ lead to 
 conformal covariants/anomalies of the four types mentioned above
(Theorem \ref{lastthm}). Applying this to explicit geometric
conformally covariant operators such as the Dirac,
 Paneitz and Peterson
operators (see \S2.2) yields conformal
anomalies and covariants, including
a new example associated to
the heat kernel asymptotics of conformally covariant
 pseudodifferential operators.  The Laurent approach provides a natural
 hierarchy among these invariants/covariants: the most divergent term in the
 Laurent expansion is a conformal invariant; 
 if this global invariant vanishes in a particular case, 
 then the new ``most divergent''  term,
if it is of the form $\int_M {\mathcal I}(g, x) \dvol_g(x)$
tends to give rise to a local conformal anomaly proportional to ${\mathcal I
}(g,x)$.


\section{Regularized  traces}

In this section, we recall known results on
regularized traces and the Wodzicki residue, and give some extensions to
families of operators.

\subsection{Preliminaries}

Let $E\to M$ be a hermitian vector bundle over a closed Riemannian
$n$-manifold $M$, and let $Cl(M, E)$ denote the algebra of classical
pseudodifferential operators ($\pdo$s)
acting on smooth sections of $E$. $S^*M\subset
T^*M $
denotes the unit cosphere bundle, and ${\rm tr}_x$
denotes the trace on the fiber $E_x$ of $E$ over $x\in M$.

\begin{defn}
A positive order
 elliptic operator $Q\in Cl(M, E)$ is {\it
 admissible} if there is an angle with vertex $0$ which contains the spectrum
 of the leading symbol $\sigma_L(Q)$ of $Q$. A choice of a half line
 $L_\theta=\{re^{i\theta}, r>0\}$ which does not intersect the spectrum of $Q$
 (which is discrete since $M$ is compact)
 is a {\it spectral cut} for Q, and $\theta$ is an {\it Agmon angle.}
An admissible operator is also
called a {\it weight}. $Ell^{ adm}_{>0}(M, E)$
 (resp. $Ell^{*, adm}_{>0}(M, E)$) is the class of 
 admissible (resp. invertible admissible)
elliptic operators of positive order in $\Cl(M, E)$.
\end{defn}

Examples of admissible elliptic operators are classical $\pdo$s
with positive leading symbol such as generalized Laplacians and
formally self-adjoint elliptic classical $\pdo$s such
as Dirac operators in odd dimensions. 

An admissible invertible elliptic
operator of positive order and with spectral cut $L_\theta$ has well-defined
complex powers (Seeley \cite{Se}) defined  for ${\rm Re}(z)$ sufficiently negative
by the contour integral
$$Q^z_\theta:= \frac{i}{2\pi} \int_{C_\theta} \lambda^z 
\, (Q-\lambda I)^{-1} \, d\lambda$$
where $C_{\theta}$ is a contour encircling $L_\theta.$
One then extends the complex power $Q_\theta^z$ to any half plane
Re $z<k, k\in \N$  via the formula $ Q_\theta^k Q_\theta^{z-k}= Q_\theta^z.$
These  complex powers  clearly   depend on the choice of spectral cut.
Setting $z=0$, we get   
$$ Q_\theta^0 = I - \Pi_Q =  \frac{i}{2\pi} \int_{C_\theta}
  \, (Q-\lambda I)^{-1} \, d\lambda,$$ 
where $\Pi_Q$ is the projection onto the generalized kernel of $Q$.
The  logarithm of $Q$, which also depends on the spectral cut, is defined by
$$\log_\theta \,Q:= \frac{d}{dz}\left|_{_{_{_{_{z=0}}}}}\right. Q_\theta^z.$$ 
This dependence  will be 
omitted  from  the notation from now on.

\subsection{The Wodzicki residue }\label{wodsubsection}

Let $A\in Cl(M,E)$ have order $\alpha$ and
symbol $\sigma(A)(x, \xi)\sim\sum_{j=0}^\infty
 \psi(\xi)\,\sigma_{\alpha-j}(A)(x, \xi), $ where $\sigma_{\alpha-j}$ is the 
positively
 homogeneous component
of order $a-j$ and $\psi$ is a smooth cut-off function which
is one outside a ball around
 $0$ and vanishes on a smaller such ball.  Let  $dx = dx^1\wedge\ldots\wedge
dx^n $ be the locally defined coordinate form on $M$, and let
$d\xi$ be
the volume form on $T^*M$ (or the restriction of $d\xi$ to
the unit cosphere bundle $S^*M\subset T^*M$ or to the unit cosphere
 $S_x^*M$ at a fixed $x\in M$).  Then
  $${\rm res}_x(A)dx:= \left(
\int_{S_x^*M} {\rm tr}_x \sigma_{-n}(A)(x,\xi) \,d\xi\right)dx,$$
is (nontrivially) a global top degree form on $M$ whose integral
$${\rm res}(A):= \frac{1}{(2\pi)^n}\,\int_M {\rm res}_x(A)\, d x$$
is the {\it Wodzicki residue \cite{W} of $A$}  
(see Kassel \cite{K}, Lesch \cite{L} for a review and further development).


The Wodzicki residue has several striking properties.  From its definition,
the Wodzicki residue vanishes on differential operators and operators of
nonintegral order, but it is nonzero in general.
The Wodzicki residue is local, in
that it is integral over $M$ of a density which is computed pointwise
from a homogeneous component of the
symbol.  Most importantly, the Wodzicki residue is cyclic on $CL(M,E)$ in the
following sense:
$$
{\rm res}([A\,, B])=0,\quad\forall A, B\in Cl(M,  E).$$

 The Wodzicki residue extends to logarithms of admissible elliptic operators
 $Q$ by
\begin{eqnarray*}
{\rm res}( \log Q) &:=& \frac{1}{(2\pi)^n}\int_M {\rm res}_x(\log
Q) dx\\
&:=& \frac{1}{(2\pi)^n} \int_{S^*M} {\rm tr}_x\sigma_{-n}( \log
Q)(x,\xi) d\xi dx
\end{eqnarray*}
(Okikiolu \cite{O1}).
More generally, given $A\in \Cl(M, E)$, if
$${\rm
res}_x(A\, \log Q)\, dx:= \left(\int_{S_x^*M} {\rm tr}_x\sigma_{-n}(A\,
\log Q)(x,\xi) d\xi\right)\, dx$$ 
defines a global form on
$M$, we can integrate it over $M$ to define
 $$
{\rm res}(A\, \log Q):= \frac{1}{(2\pi)^n} \int_{S^*M} 
{\rm tr}_x\sigma_{-n}( A\,\log Q)(x,\xi) d\xi dx.$$
This  holds in particular if   $A$ is a differential 
operator \cite[Thm.~2.5]{PS}.

 The cyclicity of the Wodzicki residue partially extends to
 logarithmic operators. The Wodzicki residue vanishes on
 brackets of the type 
 $[A, B\, \log_\theta Q]$ where $A, B\in \Cl(M, E)$, $Q\in Ell^{*,
 adm}_{0}(M, E)$, and $[A, B]$ is a differential operator \cite{O1,PS}
(Thm.~4.9).

\subsection{The canonical trace}

By a procedure well known to physicists and mathematicians (see
Paycha \cite{P} for a review),
a classical symbol $\sigma$ on $\R^n$,
has a cut-off integral in momentum space
 $\{\xi\}$. 
To set the notation, let $\psi$ be the cutoff function of
\S\ref{wodsubsection}, and set
$$\sigma_{(N)}(x, \xi):= \sigma(x,\xi)- \sum_{j=0}^N \psi(\xi) \, 
\sigma_{\alpha(z)-j}(x, \xi).$$

\begin{prop} \label{prop:finiteparts}Let $\sigma$ be a classical symbol on an
 open subset $U\subset \R^n$ of order $\alpha$.  For $x\in U$, let 
$B_x^*(R)\subset T_x^*U$ be the ball
 of radius $R$ centered at $0$. As
 $R\to \infty$,
\begin{equation}\label{3o} 
\int_{B_x^*(R)} \tr_x \sigma(x, \xi)\, d\xi\sim
\sum_{j=0, \alpha-j+n\neq 0}^\infty a_{\sigma,j} (x)\, R^{\alpha-j+n}+ 
b_\sigma(x)\,
 \log R+ c_\sigma(x),\end{equation}
with
$$a_{\sigma, j}(x)= 
{\alpha+n-j}; \quad b_\sigma(x)=  \int_{S_x^*U}\tr_x \sigma_{-n}(x,\xi)\, d\xi$$
and with  finite part/cut-off integral
 \begin{eqnarray}\label{eq:cutoffint}
c_\sigma(x) &:=& \cutoffinteq_{T_x^*U} \tr_x \sigma(x, \xi)\, d\xi\nonumber \\
&:=& {\rm fp}_{R\to \infty}
 \int_{B_x^*(R)}\tr_x \sigma(x, \xi) \, d\xi\nonumber\\
&=& \int_{T_x^*U} \tr_x \sigma_{(N)}(x, \xi)\, d\xi
+ \sum_{j=0}^N \int_{B_x^*(1)} \psi(\xi) \, 
\tr_x \sigma_{\alpha-j} (x, \xi) \, d\xi\nonumber\\
&&\qquad - \sum_{j=0, \alpha-j+n\neq 0}^\infty a_{\sigma, j}(x).
\end{eqnarray}
The finite part is independent of reparametrization of $R$ 
provided $
b_\sigma(x)$
vanishes. 
\end{prop}

Whenever $\alpha$ is nonintegral, 
via a partition of unity on $M$ one can patch the local cut-off
integrals $ \cutoffint_{T_x^*U} \tr_x \sigma_A(x, \xi)\, d\xi$  into a cut-off
integral 
\begin{equation}\label{cutoffpatch}
 \omega_{KV}(A)(x)=\cutoffinteq_{T_x^*M} \tr_x \sigma_A(x, \xi)\, d\xi
\end{equation}
on
$T_x^*M$ and then integrate over $M$ to get the Kontsevich-Vishik
{\it canonical trace} \cite{KV}
\begin{equation}\label{canon}
{\rm TR}(A):=  \frac{1}{(2\pi)^n}\int_M \omega_{KV}(A)(x)\, dx=  
\frac{1}{(2\pi)^n}\int_M dx\, \cutoffinteq_{T_x^*M}
{\rm tr}_x\sigma(A)(x, \xi)\, d\xi.  
\end{equation}
We  consider holomorphic families of classical symbols \cite{KV}. 
\begin{defn}
A  family of complex valued classical symbols $z\mapsto \sigma(z)$ 
on an open subset $U$ of $\R^n$ is  
{\it holomorphic}
 on a subset $W\subset \C$  if:
 \begin{enumerate}
 \item[1.] The order $\alpha(z)$ of $\sigma(z)$
is holomorphic\footnote{i.e. differentiable in $z$};  in $z\in W$;
\item[2.]  For any nonnegative integer $j$, 
the   map $(z,x,\xi)\mapsto \sigma(z)_{ \alpha(z) -j}(x,\xi) $
is holomorphic in $z$ and 
the map
$z\mapsto \left(\sigma(z)\right)_{ \alpha(z) -j} $ is a continuous
map from $W$ to  $\Ci(T^*U )$ in the standard topology on $\Ci(T^*U )$.
\item[3.]   For $N \gg 0$, the truncated kernel
  $$K(z)^{(N)}(x, y):= \int_{T_x^*U} e^{i \xi\cdot (x-y)} 
  \sigma(z)_{(N)}(x, \xi)d\xi$$ 
defines  a holomorphic 
 map $W\to C^{k(N)}(U\times U),
z\mapsto K(z)^{(N)}$  
for some $k(N)$ with  $\lim_{N\to \infty} k(N)=\infty$.
\end{enumerate}
A family $A(z)\in Cl(M,E)$ of classical $\pdo$s is holomorphic for
$z\in W\subset \C$ 
if it is defined in any local trivialization by a holomorphic
family of classical symbols $\sigma_{A(z)}$.
\end{defn}
The cut-off integral $\cutoffint_{T_x^*U}\tr_x \sigma(z)(x, \xi) d\xi$ is defined
whenever the order $\alpha(z)$ is nonintegral. 
The following 
extends results of Kontsevich-Vishik on the explicit Laurent expansions of
holomorphic families \cite[Thm.~2.4]{PS}:

\begin{prop} \label{prop:Laurentexp}
Let $\sigma(z)$ be a holomorphic family of classical symbols on an
open set $U\subset \R^n$ of linear order $\alpha(z)= \alpha'(0)\,
z+ \alpha(0)$ with $\alpha'(0)\neq 0$.
Then the map $z\mapsto \cutoffint_{T_x^*U} \tr_x \sigma(z)(x, \xi) d\xi$ is
meromorphic with Laurent expansion at $z=0$ given by
\begin{eqnarray*}\label{four}
\cutoffinteq_{T_x^*U}\tr_x 
\sigma(z)(x,
\xi) d\xi&=&\left(-
\frac{1}{
\alpha^\prime( 0 )}
\int_{S_x^*U} \tr_x  \sigma(0)
_{-n}(x, \xi)d\xi\right)\cdot \frac{1}{ z}\nonumber\\
&&\qquad  + \sum_{k=0}^K\frac{z^k}{k!}\left(\cutoffinteq_{T_x^*U}
\tr_x \sigma^{(k)}(0)(x,\xi)d\xi \right. \\
&&\qquad - \left .  
\frac{1}{\alpha^\prime(0)(k+1)}\int_{S_x^*U} \tr_x  \sigma^{(k+1)}(0)
_{-n}(x, \xi)d\xi\right)\nonumber\\
&&\qquad + {\rm O}(z^K),\nonumber\\
\end{eqnarray*}
for $K\geq 0.$ 
Applying this to the symbols $\sigma_{A(z)}$ of a holomorphic family
$A(z)$ of classical $\pdo$s,
taking the
fibrewise trace and replacing  $U$ by $M$ via a partition of the
unity
provides an analogous formula 
for the first $k+1$ terms of
the Laurent expansion around $0$ of 
 $\omega_{KV}(A(z))(x)$ defined by (\ref{cutoffpatch}) with $A$
replaced by $A(z)$ and hence, after integration over $M$, of the
canonical trace $\TR(A(z))$.
\end{prop}

\begin{rem} 
(i)
Even though $\sigma(z)$ is a classical symbol, 
$\sigma^{(k)}(0)$ need not be \cite{PS}.  

(ii)  
 Since
$\omega_{KV}(A(z))(x) dx= \left(\cutoffint_{T_x^*M}\tr_x \sigma_{A(z)} (x,
\xi) \, d\xi\right) dx$ defines a global form, the coefficient of $z^k$ in
the Laurent expansion of  $\omega_{KV}(A(z))(x) dx$ 
also gives rise to a
global form.

(iii) 
For a classical $\pdo$ $A$ of order $\alpha$, 
the operator $A(z)= AQ_\theta^{-z}$ defines a
holomorphic family of classical $\pdo$s
of order $\alpha(z)=- q
z+\alpha$.  From Proposition \ref{prop:Laurentexp}, 
we recover the well known result relating
the Wodzicki residue to a complex residue:
\begin{equation}\label{starstar}
{\rm Res}_{z=0}\omega_{KV}(AQ^{ -z})(x)= -\frac{1 }{q} \res_x(A),
\end{equation}
which after integration over $M$ yields
$${\rm Res}_{z=0}{\rm  TR}(AQ^{-z}) = -\frac{1}{q} {\rm res}(A).$$
If $A$ is a differential operator, 
$\cutoffint_{T_x^*M}\tr_x \sigma^{}(0)(x,\xi)d\xi=
\cutoffint_{T_x^*M}\tr_x\sigma_{A(0)}(x,\xi)d\xi $ vanishes.
Therefore
$\left(\int_{S_x^*M}\tr_x \left(\sigma^\prime(0)\right)_{-n}(x,
\xi)\right)d\xi dx$ defines a global form, whose integral is
$-{\rm
res}\left(A\,\log_\theta Q\right)$ \cite{PS}.    
\end{rem} 

\subsection{Weighted traces}

For a weight $Q$ with spectral cut and a nonnegative integer
$k$, set
$${\mathcal A}^k(M, E):= \{\sum_{j=0}^k A_j\log^jQ, \, A_j\in Cl(M, E),
\quad 0\leq j\leq k\}.$$ 
Operators in ${\mathcal A}^k(M, E) $ coincide with Lesch's
log-polyhomogeneous operators \cite{L}.\\
  ${\mathcal
A}^k(M, E) $ is in fact independent of the choice of $Q$  (Ducourtioux
 \cite{D,PS}) and coincides with the class $Cl^{*, k}(M, E)$ of Lesch
 \cite{L}.  Note that ${\mathcal A}^0(M, E)= Cl(M, E)$.
The order of $A_j\, \log^j Q$ is defined to be
the order of $A_j$. 

Cut-off
integrals extend \cite{L} to symbols of operators in ${\mathcal A}^k(M,
E)$, once (\ref{3o}) is extended to include the terms $d_{\sigma,j}
\log^j R, j=1,
\ldots, k+1$.  As for classical operators, for a noninteger
order $A\in {\mathcal A}^k$, $\omega_{KV}(A)(x)\, dx:=
\cutoffint_{T_x^*M} {\rm tr}_x(\sigma(A)(x, \xi) \, d\xi$ defines a global
form, and one can define the  canonical trace ${\rm TR}(A)$ by
(\ref{canon}).
The linear functional TR is  cyclic:
\begin{equation} \label{eq:TRcyclic}
{\rm TR}([A, B])
=0,\quad \forall A\in {\mathcal A}^k(M, E), B\in {\mathcal A}^j(M, E),
\quad  {\rm ord}(A)+ {\rm ord} (B)\notin \Z.
\end{equation}


Weighted traces are defined by the finite part in the Laurent expansion of the
canonical trace of a holomorphic family; this is in contrast to
the Wodzicki residue, which occurs as the residue in the Laurent expansion.
\begin{defn}
For $A\in {\mathcal A}^k(M, E)$, the {\it $Q$-weighted trace} of $A$ is
\begin{eqnarray*}
{\rm tr}^Q(A)&:=& {\rm fp}_{z=0} 
{\rm TR}(AQ^{-z})+{\rm tr}(A\,\Pi_Q)\\
&:=& \lim_{z\to 0} \left( {\rm TR}(AQ^{-z})-\sum_{j=0}^{k} 
\frac{a_{j+1}}
{z^{j+1}}\right)+{\rm tr}(A\,\Pi_Q),
\end{eqnarray*}
where $ a_{j+1} $ is the residue of $ {\rm TR}(AQ^{-z})$ of order $j+1$. 
\end{defn}

The existence of the Laurent expansion is known \cite{L}.
As usual, this definition depends on a choice of spectral cut
    for $Q$.  For $A\in \Cl(M,E)$, the 
weighted trace can also be defined by the finite part of
    ${\rm tr}(AQ^{-z})$, where tr is the ordinary operator trace.  Indeed, for
    ${\rm Re}(z)\gg 0$, $AQ^{-z}$ is trace-class, in which case ${\rm
    TR}(AQ^{-z}) = {\rm tr}(AQ^{-z}).$  The known meromorphic
    continuation of the right hand side (Grubb and Seeley \cite{GS})
gives the equivalence of the two
    definitions.  We prefer our current definition of the weighted trace, 
since ${\rm TR}(AQ^{-z})$ is well defined outside a countable set of poles,
and hence does not require a meromorphic continuation.


Weighted traces do not have the local properties of the Wodzicki
residue in general.  For example,
a formally self-adjoint, positive order,
invertible elliptic operator $A\in \Cl(M, E)$ is admissible
with Agmon angle
$\theta=\frac{\pi}{2}$, as is its modulus $\vert
A\vert:=\sqrt{A^*A},$ which has positive leading symbol. Then $A\,\vert
A\vert^{-1}\in\Cl(M, E)$, and we can set 
$$\eta_A (z):= {\rm
TR}\left (A\,\vert A\vert^{-z-1}\right).$$ The $\eta$-invariant of $A$ is
given by its finite part:
\begin{equation}\label{star1}
\eta_A(0):= {\rm tr}^{\vert A\vert} (A\vert A\vert^{-1}),\end{equation}
which is not local in general.

\begin{rem}
The map $z\mapsto \eta_A (z)$ is holomorphic at $z=0$ since the
Wodzicki residue of a $\pdo$ projection such as ${\rm res}(A\,\vert
A\vert^{-1})$
vanishes.  It follows 
that $\eta_A(0)
= {\rm tr}^{ A} (A\vert A\vert^{-1})$, {\it i.e.}~the
$\eta$-invariant can be defined using the easier $A$ 
as a weight (Cardona, Ducourtioux and Paycha \cite[Prop.~1]{CDP}.
\end{rem}

For differential operators $A$,  ${\rm res} (A\,
 \log Q)$
is well defined \cite[Thm.~3.7]{PS}, and 
\begin{equation}\label{eq:trQres}
{\rm tr}^Q(A) = -\frac{1}{q} {\rm res} (A\, \log Q)
= -\frac{1}{q } \int_{M}{\rm res}_x (A\, \log Q) dx.
 \end{equation} 
In this case, ${\rm tr}_\theta^Q(A)$ has a partial
locality as an integral of 
$\sigma_{-n}(A\,\log Q)$. In
particular, for $A=I$ we have 
\begin{equation}\label{quick} 
{\rm tr}^Q(I)= -\frac{1}{q} {\rm res} (\log Q),\end{equation}
an
expression related to the exotic determinant ${\rm
det}_{res}(Q)= e^{{\rm res}(\log Q)}$ \cite{Sc} (and references
therein).  In turn,
\begin{equation}\label{star2}
{\rm tr}_\theta^Q(I)-{\rm tr}(\Pi_Q) = \zeta_Q(0),\end{equation}
 where the zeta function
is given by the usual meromorphic continuation of
$$\zeta_{Q}(z)= {\rm  TR} (Q_\theta^{-z}) = {\rm tr}(Q^{-z}),$$
which is well defined 
for ${\rm Re}(z) > \frac{n}{q}$ ($n = {\rm dim}(M),\ q = {\rm ord}(Q)>0$).

Since the Wodzicki residue vanishes for differential operators $Q$, $\zeta_{Q
}(z)$ is holomorphic at $z=0$, and an easy computation yields
\begin{equation}\label{star3}
 \zeta_{Q}^\prime(0)= -{\rm tr}^Q(\log Q)\end{equation}
for an invertible weight $Q$.\\

In summary, the key spectral invariants $\eta_A(0), \zeta_Q(0), \zeta'_Q(0)$
all occur as weighted traces.

The following proposition 
will be used in \S2.

\begin{prop}\label{prop:jterm} Let $A\in \Cl(M,
 E)$ and let $Q$ be an invertible weight. We have the Laurent expansion
$${\rm TR}(A\, Q^{-z})= \frac{{\rm res}(A)}{q\, z} +
\sum_{j=0}^J \frac{(-1)^j}{j!} \, 
{\rm tr}^Q(A\, \log^j Q)\, z^j+ {\rm o}(z^J).$$
\end{prop}

\noindent{\bf Proof:} By Remark 1, the map $z\mapsto
{\rm TR}(A\, Q^{-z})$ is meromorphic with a simple pole at $z=0$ with
residue  $ \frac{{\rm res}(A)}{q}$, so 
$${\rm TR}(A\, Q^{-z})= \frac{{\rm res}(A)}{q\, z} +\sum_{j=0}^J a_j(A, Q)\,
z^j+ {\rm o}(z^J).$$ \\
Since Laurent expansions can be differentiated term by term away from their
poles, we obtain
\begin{eqnarray*}
{\rm tr}^Q(A\, \log^j Q) &=& {\rm fp}_{z=0}{\rm TR}(A\,\log^j Q\, Q^{-z})
=
(-1)^j\, {\rm fp}_{z=0} \left(\partial^j_z{\rm TR}(A\,  Q^{-z})\right)\\
&=& (-1)^j\, j!\, a_{j}(A, Q). \hfill \Box
\end{eqnarray*}

\begin{rem} 
$z\mapsto {\rm TR}(A\,\log^j Q\,
Q^{-z})$ has a Laurent expansion \cite{L} with poles of order at most $j+1$:
$${\rm TR}(A\,\log^j Q\,  Q^{-z})= \sum_{l=1}^{j+1}
\frac{ b_{j,l}(A, Q) }{z^{l}}+ \sum_{i=0}^k a_{j, i}(A, Q)\, z^i+ {\rm o}(z^k).$$
The $a$ and $b$ coefficients are related.
For example, the identity $\partial_z \TR(A\log^j Q
Q^{-z})= -\TR(A\log Q^{j+1} Q^{-z})$ implies
$$a_{j+1,i}=-(i+1)a_{j, i+1}(A, Q), b_{j+1, l+1}(A, Q)=l\, b_{j, l}(A,Q).$$
\end{rem}

\subsection{Differentiable families of canonical traces}

The definition of a $C^k$ differentiable family of classical symbols is
completely analogous to the the holomorphic definition.  Namely, the
one-parameter family of symbols
$t\mapsto \sigma_t, t\in \R$, with $\sigma_t$
 defined on an open set $U\subset \R^n$, is $C^k$  for a fixed
$k\in \Z^+$  if (i)
the order $\alpha_t$ of $\sigma_t$ is $C^k$ in $t$, (ii) each homogeneous
component $\sigma_{t, \alpha_t -j} (x, \xi)$ is $C^k$ in $t$, 
(iii) for 
$N\gg 0$ and $K_t^{(N)}(x, y)$ the truncated kernel,   
the map $U\to C^{K(N)},\ 
t\mapsto K_t^{(N)}(x,y)$ is
$C^k$ for some $K(N)$ with 
$\lim_{N\to \infty} K(N)=\infty$.
A family $t\mapsto A_t$ of classical $\pdo$s is  $C^k$  if it is defined  
in any local trivialization by  $C^k$ family of symbols.

\begin{rem} \label{rmkfour} {\rm 
By (iii), for a 
 $C^k$ family $t\mapsto \sigma_t$, the map $t\mapsto
 \left(\sigma_t\right)_{(N)}$ is also  $C^k$
 in $\Ci(T^*U)$, and then by (ii), this differentiability
holds for $t\mapsto \sigma_t$.  As a consequence, for a $C^k$
family
 $t\mapsto \sigma_t$ and for any compact set $K\subset
 T^*U$, $t\mapsto \Vert\partial_t^k\sigma_t \Vert_K:= {\rm sup}_{(x, \xi)\in
 K}\vert\partial_t^k \sigma_t(x,\xi)\vert$ is continuous and hence uniformly
 bounded on any interval $[t_0-\eta, t_0+\eta]$, $\eta>0,$ as are
 the homogeneous components $\left(\sigma_t\right)_{\alpha-j}$.  Moreover,
the minus one order symbol $(\vert
 \xi \vert+1)^{N-\alpha}\left(\partial_t^k\sigma_t\right)_{(N)}(x, \xi)$ 
is  bounded on $T_x^*U$ and gives rise to a
 continuous map 
\begin{eqnarray*} t &\mapsto& (\vert \xi
 \vert+1)^{N-\alpha}\Vert\left(\partial_t^k\sigma_t\right)_
{(N)}\Vert_{T_x^*U}\\
&&\qquad :=
 {\rm sup}_{\xi \in T_x^*U} \left( (\vert \xi
 \vert+1)^{N-\alpha}\vert\left(\partial_t^k\sigma_t\right)_{(N)}(x,\xi)
\vert\right). 
\end{eqnarray*}
Hence $(\vert \xi \vert+1)^{N-\alpha} \vert
 \left(\partial_t^k\sigma_t\right)_{(N)}\vert$ is
 uniformly bounded above on $ [t_0-\eta, t_0+\eta]$ by a constant
 $C_{t_0, \eta, (N)}$. Therefore, for fixed $x$ the map $\xi\mapsto \vert
 \left(\partial_t^k\sigma_t\right)_{(N)}(x,\xi)\vert$ 
is bounded above by a map
 $\xi \mapsto C_{t_0, \eta, (N)}\, (\vert \xi \vert+1)^{\alpha-N},$ which lies
 in $L^1(T_x^*U)$ for $N \gg 0.$  }
\end{rem}

This remark implies that the cut-off integral and the canonical trace
commute with differentiation as long as the symbols and operators have
constant noninteger order.

\begin{thm} \label{thm:TRdiff}
\begin{enumerate}
\item[1.] Let $t\mapsto \sigma_t$ be a $C^1$ family of symbols on $U$ with
constant noninteger order $\alpha$. Then
$$     
\frac{d}{dt}\cutoffinteq_{T_x^*U} \tr_x \sigma_t(x, \xi) \, d\xi 
=\cutoffinteq_{T_x^*U} \tr_x \dot \sigma_t(x, \xi),
$$
where $\dot\sigma_t = \frac{d}{dt}\sigma_t.$
\item[2.] Let $t\mapsto A_t\in\Cl(M, E)$ be a $C^1$ family 
of constant noninteger order operators. Then
\begin{equation} \label{eq:TRdiff}
\frac{d}{dt}{\rm TR}(A_t)={\rm TR}(\dot A_t).
\end{equation}
\item[3.] 
Assume that for fixed $t$, $z\mapsto \sigma_t(z)$ is a holomorphic
family of classical symbols on $U$ parametrized by $z\in W\subset \C$ with
holomorphic order $\alpha(z)$ independent of $t$ and that $t\mapsto
\sigma_t(z)$ is a $C^1$ family for fixed $z\in W$. Then $z\mapsto
\cutoffint_{T_x^*U}\tr_x \sigma_t(z)(x, \xi) \, d\xi $ and $z\mapsto
\cutoffint_{T_x^*U}\tr_x 
\dot \sigma_t(z)(x, \xi) \, d\xi$ are meromorphic in $z$,
and the Laurent expansion of $\cutoffint_{T_x^*U}\tr_x 
\dot \sigma_t(z)(x, \xi)
\, d\xi$ around $z=0$ is obtained by term by term $t$-differentiation
of the Laurent expansion of
$\cutoffint_{T_x^*U}\tr_x  \sigma_t(z)(x, \xi) \, d\xi$.
\item[4.] 
 Assume that for fixed $t$, $z\mapsto A_t(z)\in Cl(M, E)$ defines a
holomorphic family on $W\subset \C$ with holomorphic order $\alpha(z)$
independent of $t$, and assume that $t\mapsto
\partial_z^k|_{z=0}A_t(z)$ is a $C^1$ family for $k\in \Z^{\geq 0}.$
Then $z\mapsto {\rm TR}(A_t(z))$ and $z\mapsto {\rm TR}(\dot A_t(z))$
are meromorphic in $z$, and the Laurent expansion of ${\rm TR}(\dot
A_t(z))$ around $z=0$ is obtained by 
term by term $t$-differentiation
of the Laurent expansion of
${\rm TR}( A_t(z))$.
\end{enumerate}
\end{thm}

{\bf Proof:} 
\noindent 1. 
Once we justify pushing the derivative past the integral, by 
(\ref{eq:cutoffint}) (and noting that we may choose
$N$ independent of $t$ by our assumption on $\alpha(z)$), we have
\begin{eqnarray*}
\lefteqn{
\frac{d}{dt} \cutoffinteq_{T_x^*U}\tr_x  \sigma_t(x, \xi)\, d\xi}\\
 &=&
\frac{d}{dt}\int_{T_x^*U}\tr_x \left( \sigma_t\right)_{(N)}(x, \xi)\, d\xi+
\frac{d}{dt}\sum_{j=0}^N \int_{B_x^*U} \psi(\xi)
\tr_x \left(\sigma_t\right)_{\alpha-j} (x, \xi) \, d\xi\\ 
&&\qquad +\frac{d}{dt}
\sum_{j=0, \alpha-j+n\neq 0}^\infty
\frac{\int_{S_x^*U}\tr_x \left(\sigma_t\right)_{\alpha-j}(x,\xi)\,
d_S\xi}{\alpha+n-j}\\ 
&=& \int_{T_x^*U}\frac{d}{dt}\tr_x 
\left(\sigma_t\right)_{(N)}(x, \xi)\, d\xi+ \sum_{j=0}^N \int_{B_x^*U}
\psi(\xi)
 \frac{d}{dt}\tr_x
\left(\sigma_t\right)_{\alpha-j} (x, \xi) \, d\xi\\
&&\qquad + \sum_{j=0, \alpha-j+n\neq 0}^\infty
\frac{\int_{S_x^*U}\frac{d}{dt}\tr_x \left(\sigma_t\right)_{\alpha-j}(x,\xi)\,
d_S\xi}{\alpha+n-j}\\ 
&=& \cutoffinteq_{T_x^*U}\tr_x  \dot \sigma_t(x, \xi)\,
d\xi,
\end{eqnarray*}
Recall that
for $\xi\in A\subset \R^n$, $t_0\in \R$ and  $\epsilon >0$, if
 $\vert t-t_0\vert \leq \epsilon$  implies $\vert
\frac{d}{dt}f(t, \xi)\vert \leq g(\xi) $ with $g\in L^1(A)$, then
$\frac{d}{dt}\vert_{_{t=t_0}} \int_A f(t, \xi) \, d\xi = \int_A
\frac{d}{dt}\vert_{_{t=t_0}} f(t, \xi) \, d\xi.$ This applies to
the
compact subsets $A=B_x^*U$ and $A= S_x^*U$ and
$f(t,\xi)=\tr_x \left( \sigma_t\right)_{\alpha-j}(x, \xi)$ and to 
$A=T_x^*U$
and $f(t, \xi)=
\tr_x \left(\sigma_t\right)_{(N)}(x, \xi)$ (where the 
the required uniform estimates follow from 
Remark \ref{rmkfour} with
$k=1$).

\noindent 2. 
By 1, 
\begin{eqnarray}\label{eq:Trdiffproof}
\frac{d}{dt} \cutoffinteq_{T_x^*M}
{\rm tr}_x\sigma(A_t)(x, \xi)\, d\xi &=& \cutoffinteq_{T_x^*M}
{\rm tr}_x\dot \sigma( A_t)(x, \xi)\, d\xi\nonumber\\
&=& \cutoffinteq _{T_x^*M}
{\rm tr}_x\sigma(\dot A_t)(x, \xi)\, d\xi.
\end{eqnarray}
 Since $A_t$ has constant order $\alpha$, $\dot A_t$ has order $\alpha$ modulo
integers. Therefore $\dot A_t$ has noninteger order, so 
$ \left(\cutoffint_{T_x^*M} {\rm tr}_x\sigma(\dot A_t)(x, \xi)\,
d\xi\right)\, dx$ is a global form on $M$ and ${\rm TR}(\dot A)=
\frac{1}{(2\pi)^n}\int_M dx\, \cutoffint_{T_x^*M} {\rm tr}_x\sigma(\dot A)(x,
\xi)\, d\xi $ is well defined. 
 Since\\
 ${\rm TR}(A)=  \frac{1}{(2\pi)^n}\int_M dx\, \cutoffint_{T_x^*M}
{\rm tr}_x\sigma(A)(x, \xi)\, d\xi,  $  
integrating  (\ref{eq:Trdiffproof}) over
$M$  yields  (\ref{eq:TRdiff}).

\noindent 4.  We now prove 4, leaving the similar proof of 3 to the
reader. 
If 
$A_t(z)$ has noninteger order, by 2
$$\frac{d}{dt} {\rm TR}(A_t(z))= {\rm TR}(\dot A_t(z))$$ except at the poles,
so this is an equality of meromorphic functions. By Proposition
\ref{prop:Laurentexp}, the coefficients of the Laurent expansion on either
side can be expressed in terms of the cut-off integral of
$\tr_x \sigma\left(\partial_z^k A_t(z)\vert_{{z=0}}\right)$
(resp. $\tr_x \sigma(\partial_z^k \dot A_t(z)\vert_{{z=0}})$) on
$T_x^*U$ and ordinary integrals over compact sets of the $-n$ component
of $\tr_x \partial_z^j \sigma(A_t(z))\vert_{_{z=0}}$
(resp. $\tr_x \partial_z^j \sigma(\dot A_t(z))\vert_{_{z=0}}$),
$j\in \Z^{\geq 0}$. 
As above
$t\mapsto \left(\partial_z^k A_t(z)\right)\vert_{_{z=0}} $ is $C^1$, 
so we can push the 
derivative past the integral as desired.
{\hfill $\Box$}\\

Let $h:W\subset \C\to \R$ be a $C^1$ map such that
$$h(A):=\frac{i}{2\pi} \int_{C_\theta} h(\lambda)\, (A-\lambda)^{-1}\,
d\lambda; \quad h^\prime(A):=\frac{i}{2\pi} \int_{C_\theta}
h^\prime(\lambda)\, (A-\lambda)^{-1}\, d\lambda
 $$ 
takes  any weight $A$  to $h(A),
h^\prime(A)\in \Cl(M,E)$. Here $\theta$ is an Agmon angle 
for $A$ and
$C_\theta$ is the associated contour
(where we
assume  $C_\theta\subset W$). Examples of such maps are
$$h(z)=\frac{z}{\vert z\vert},\quad W=\R^*;\quad
h(z)= z_\theta^c,\quad W=\C;\quad  h(z)=\frac{z}{\vert z\vert},\quad W=\R^*. 
 $$ 
for fixed $c\in \R.$

\begin{prop}
Let $t\mapsto A_t$ be a differentiable family of weights of constant
noninteger order
and with common Agmon angle.
Then
\begin{equation}\label{eq:TRh}
\frac{d}{dt} {\rm TR}\left( h(A_t) \, A_t^{-z}\right)= {\rm TR}\left(
h^\prime(A_t) \, \dot A_t\, A_t^{-z}\right)-z \, {\rm TR}\left( h(A_t) \, \dot
A_t\, A_t^{-z-1}\right).\end{equation} 
This is equivalent to the
following set of equations:
\begin{eqnarray}\label{eq:TRhresterm}
 \frac{d}{dt} {\rm res}\left(h(A_t) \right)
&=&  {\rm res}\left(  h^\prime(A_t)\, \dot A_t \right),\\
\label{eq:TRh0term}
  \frac{d}{dt} {\rm tr}^{A_t}\left(h(A_t) \right)
&=&  {\rm tr}^{A_t}\left(  h^\prime(A_t)\, \dot A_t \right)-
\frac{1}{q} {\rm res}\left( h(A_t) \, \dot A_t\, A_t^{-1}\right),\\
\label{eq:TRhjterm}
  \frac{d}{dt} {\rm tr}^{A_t}\left(h(A_t) \,\log^j A_t\right)
&= & {\rm tr}^{A_t}\left(  h^\prime(A_t)\, \dot A_t \,\log^j
  A_t\right)\nonumber \\
&&\qquad +j \, {\rm tr}^{A_t}\left( h(A_t) \, \dot A_t\, A_t^{-1}\,\log^{j-1}
  A_t\right)\end{eqnarray}
for $j\in \Z^+.$
\end{prop}

\noindent {\bf Proof:} Applying  Theorem \ref{thm:TRdiff}
gives the following equalities of meromorphic functions:
\begin{eqnarray}
\lefteqn{
\frac{d}{dt} {\rm TR}\left( h(A_t) \, A_t^{-z}\right) = {\rm
TR}\left(\frac{d}{dt}\left( h(A_t) \, A_t^{-z}\right)\right)}\\
&=& {\rm
TR}\left( \frac{d}{dt}\left( h(A_t)\right) \, A_t^{-z}\right)+ {\rm TR}\left(
h(A_t) \,\frac{d}{dt}\left( A_t^{-z}\right)\right)\nonumber\\ 
&=& -\frac{i}{2\pi} {\rm
TR}\left( \int_{C_\theta} h(\lambda)\, (A_t-\lambda)^{-1}\,\dot A_t\,
(A_t-\lambda)^{-1}\, d\lambda \, A_t^{-z}\right)\nonumber\\ 
&&\qquad - \frac{i}{2\pi} {\rm
TR}\left(h(A_t) \int_{C_\theta} \lambda^{-z}\, (A_t-\lambda)^{-1}\,\dot A_t\,
(A_t-\lambda)^{-1}\, d\lambda \right)\nonumber\\
&=& -\frac{i}{2\pi} {\rm TR}\left(
\left(\int_{C_\theta} h(\lambda)\, (A_t-\lambda)^{-2}\, d\lambda\right)\,\dot
A_t \, A_t^{-z}\right)\nonumber\\ 
&&\qquad - \frac{i}{2\pi} {\rm TR}\left(h(A_t)\left(
\int_{C_\theta} \lambda^{-z}\, (A_t-\lambda)^{-2}\, d\lambda \right)\,\dot A_t
\right)\label{cyclicity}\\  
&=&
\frac{i}{2\pi} {\rm TR}\left( \left( \int_{C_\theta} h^\prime(\lambda)\,
(A_t-\lambda)^{-1}\, d\lambda\right)\,\dot A_t \, A_t^{-z}\right)\nonumber\\
&&\qquad -\frac{i\, z}{2\pi} {\rm TR}\left(h(A_t)\left( \int_{C_\theta} 
\lambda^{-z-1}\,
(A_t-\lambda)^{-1}\, d\lambda \right)\,\dot A_t \right)\label{intparts}\\ 
&=& {\rm TR}\left(
h^\prime(A_t)\,\dot A_t \, A_t^{-z}\right)-z {\rm TR}\left(h(A_t)\,
A_t^{-z-1}\,\dot A_t\right)\nonumber\\ 
&=& {\rm TR}\left( h^\prime(A_t)\,\dot A_t \,
A_t^{-z}\right)-z {\rm TR}\left(h(A_t)\, \dot A_t\ A_t^{-z-1}
\right).\label{cyclicity2}
\end{eqnarray}

In (\ref{cyclicity}), (\ref{cyclicity2}), we use the cyclicity of TR on
noninteger order operators, and in (\ref{intparts}) we use integration by
parts.  This proves (\ref{eq:TRh}).

 By Theorem \ref{thm:TRdiff}.3, the Laurent expansion of $\frac{d}{dt}
{\rm TR}\left( h(A_t) \, A_t^{-z}\right)$ is the term by term derivative 
of the Laurent expansion of
$ {\rm TR}\left( h(A_t) \, A_t^{-z}\right)$ .  The rest of the Proposition 
then follows from identifying the coefficients in the Laurent expansions in
(\ref{eq:TRh}) and using Proposition \ref{prop:jterm}.
{\hfill $\Box$}\\


\section{Conformal invariants and anomalies}
In this part of the paper, we use canonical traces to
build functionals of conformally
covariant operators and study their conformal
properties.

\subsection{The conformal anomaly and associated two-tensor}

   Let $M$ be a Riemannian manifold and  $\Met(M)$ denote the space of
   Riemannian metrics on $M$. $\Met(M)$ is trivially
a Fr\'echet manifold as
the open cone of positive definite
   symmetric (covariant) two-tensors inside the Fr\'echet space
$$\Ci(T^*M\otimes_s T^*M):= \{h\in C^\infty(T^*M\otimes T^*M):
   h_{ab}=h_{ba}\}$$ 
of all smooth symmetric two-tensors.
The Weyl group  $W(M):= \{e^f: f\in \Ci(M)\}$ 
acts smoothly on $\Met(M)$  by
Weyl transformations 
$$W(g, f) = \bar g:=e^{2f} g,$$
and given a reference metric $g\in \Met(M)$, a  functional 
${\mathcal F}:\Met(M)\to \C$  induces a map
$$
{\mathcal F}_g= {\mathcal F}\circ W(g, \cdot):\Ci(M)\to  \C, \ \ 
f\mapsto  {\mathcal F}(e^{2f} g).$$

\begin{defn}
A functional ${\mathcal F}$ on $\Met(M)$ is {\it conformally invariant for a
 reference metric} $g$ if ${\mathcal F}_g$ is constant on a conformal class,
 i.e.
$${\mathcal F}(e^{2f} g)={\mathcal F}( g)\quad \forall f\in \Ci(M).$$ 
A
functional ${\mathcal F}$ on $\Met(M)$ is {\it conformally invariant} if it is
conformally invariant for all reference metrics.  A functional ${\mathcal F}:
\Met(M)\times M\to\C$ is called a {\it pointwise
conformal covariant of weight  w} if
$${\mathcal F}(e^{2f}g,x) = w\cdot f(x) 
{\mathcal F}(g,x)\quad \forall f\in\Ci(M),
\quad \forall x\in M.$$
\end{defn}

For conformal covariants, we always assume that ${\mathcal F}(g,x)$ is given
by a universal formula in the components of $g$ and their derivatives at $x$.

For a fixed Riemannian metric $g = (g_{ab})$, $C^\infty(M)$ has the $L^2$
metric 
$$(f,\tilde f)_g = \int_M f(x)\tilde f(x) \dvol_g(x).$$
This extends to 
the $L^2$ metric on $\Met(M)$ 
given by
\begin{equation}\label{ip}
\langle h,k\rangle_g:= 
\int_M g^{ac}(x) g^{bd}(x)\, h_{ab}(x) \, k_{cd}(x)\, \dvol_g(x),
\end{equation}
with $(g^{ab}) = (g_{ab})^{-1}.$
The $L^2$ metric induces a weak $L^2$-topology on  $\Met(M)$, and 
$L^2(T^*M\otimes_s T^*M)$, the $L^2$-closure of $\Ci(T^*M\otimes_s
T^*M)$ with respect to $\langle\ ,\ \rangle_g,$ 
is independent of the choice of
$g$ up to Hilbert space isomorphism.  The choice of a reference
metric yields the
inner product (\ref{ip}) on the tangent space
$T_g\Met(M) =\Ci(T^*M\otimes_s T^*M)$, giving the weak $L^2$ Riemannian metric
on $\Met(M)$, and forming the completion of each tangent space.  

The metric $g$ allows us to contract a two-tensor via
$${\rm tr}_g(h):= h^b_b = g^{ab} h_{ab}.$$  
The various inner products are related as follows:

\begin{lem} \label{2lem}
For $g\in\Met(M), h\in \Ci(T^*M\otimes_s T^*M)$ and 
$f \in \Ci(M)$, we have
$$\langle  h, f\, g\rangle_g= \left( {\rm tr}_g(h),f\right)_g.$$
\end{lem}
{\bf Proof:}
We have 
\begin{eqnarray*}
\langle h, f\, g\rangle_g&=& \int_M f(x)\, g^{ac}(x) g^{bd}(x)\,
h_{ab}(x) \, g_{cd}(x)\, \dvol_g(x)\\ 
&=& \int_M f(x)\,
g^{ab}(x) \, h_{ab}(x) \, \dvol_g(x)\\ 
&=& \left( {\rm tr}_g(h),f\right)_g. \hskip 2.5in \Box
\end{eqnarray*}
\medskip

A functional ${\mathcal F}: \Met(M)\to \C$ which is Fr\'echet
differentiable 
has a differential
$$ d{\mathcal F}(g): T_g\Met(M) = 
\Ci(T^*M\otimes_s T^*M)\to  \C,$$
$$d{\mathcal F}(g).h:=  \frac{d}{dt}\left|_{_{_{_{_{t=0}}}}}\right.
\frac{{\mathcal F}( g+th)-{\mathcal F}(g)}{t}.$$
For such an ${\mathcal F}$, the differentiability of
the Weyl map implies that
the composition ${\mathcal F}_g:\Ci(M)\to \C$ is differentiable at $0$ 
with differential $d{\mathcal F}_g(0):T_0\Ci(M) = \Ci(M)\to \C$.\\

\begin{defn}
The {\it conformal anomaly} for the reference metric $g$ of a differentiable
functional ${\mathcal F}$ on $\Met(M)$ is $d{\mathcal F}_g(0)$.  In physics
notation, the conformal anomaly in the direction $f\in \Ci(M)$ is
\begin{eqnarray*}\label{eq:confanomF}
\delta_f{\mathcal F}_g &:=& d{\mathcal F}_g(0).f
= d{\mathcal F}(g).2f\, g\nonumber\\
&=& \lim_{t=0}\frac{{\mathcal F}( g+2tf g)-{\mathcal F}(g)}{t}
=    \frac{d}{dt}\left|_{_{_{_{_{t=0}}}}}\right.
{\mathcal F}( e^{2tf} g). 
\end{eqnarray*}
\end{defn}

\begin{rem} 
${\mathcal F}$ is  conformally invariant if and only if
$d{\mathcal F}_g(0).f= 0$
for all $g\in \Met(M),  f\in \Ci(M).$
\end{rem}

If the differential $d{\mathcal F}(g): \Ci(T^*M\otimes_s T^*M)\to \C$
extends to a continuous functional $\overline{ d{\mathcal F}(g)}:
L^2(T^*M\otimes_s T^*M)\to \C$, then by Riesz's lemma there is a unique
two-tensor $T_g({\mathcal F})$ with
$$ \overline{ d{\mathcal F}(g)}.h=\langle h, T_g({\mathcal F})\rangle_g,
\quad \forall h\in L^2(T^*M\otimes_s T^*M).$$
$T_g({\mathcal F})$ is precisely the $L^2$ gradient of ${\mathcal F}$ at $g$.

\begin{prop}
Let ${\mathcal F}$ be a functional on $\Met(M)$ which is differentiable at the
metric $g$ and whose differential $d{\mathcal F}(g)$ extends to a
continuous functional $\overline{ d{\mathcal F}(g)}: L^2(T^*M\otimes_s T^*M)\to
\C$. Then the differential $d{\mathcal F}_g(0)$ 
also extends to a continuous functional
$\overline{d{\mathcal F}_g(0)}: L^2(M)\to \C $.  Identifying the conformal
anomaly at $g$ with a function in $L^2(M)$, we have
$$\overline{ d{\mathcal F}_g(0)}=2\, {\rm tr}_g\left( T_g({\mathcal F})\right).$$
In particular, the functional ${\mathcal F}$ is conformally invariant 
iff  $ {\rm tr}_g\left( T_g({\mathcal F})\right) = 0$ for all metrics $g$. 
\end{prop}

{\bf Proof:} The differential $d({\mathcal F}_g)_0$ extends to 
a continuous functional
because
$$ d{\mathcal F}_g (0).f=   d{\mathcal F}(g) (2f\ g) \Rightarrow
\overline{d{\mathcal F}_g (0)}.f = \overline{ d{\mathcal F}(g)}.2f\ g$$
By Lemma \ref{2lem},
$$ \overline{d{\mathcal F}_g (0)}.f=
           \overline{ d{\mathcal F}(g)}.2f\ g
= \langle  T_g({\mathcal F}) , 2f\, g\rangle_g
= 2\left( {\rm tr}_g(T_g({\mathcal F})),f\right)_g,
$$
as desired. \hfill $\Box$\\

\begin{defn} Under the assumptions of the Proposition, the function
$$x\mapsto \delta_x {\mathcal F}_g:=2 {\rm tr}_g\left( T_g({\mathcal F})\right)(x)$$
is called the {\it local anomaly} of the functional ${\mathcal F}$ at the
 reference
metric $g$.
\end{defn}

\begin{exm}
 In field theory, for a classical action ${\mathcal A}$ 
on a configuration space ${\rm Conf}(M)$ with respect to a background metric
 $g$, 
${\mathcal A}:
\phi \mapsto  {\mathcal A}(g)( \phi),$
where ${\rm Conf}(M)$ is usually a space of tensors on
$M$, the associated two-tensor $T_g({\mathcal A})$ is called the {\it
stress-energy momentum tensor}.  In the path integral approach to quantum
field theory, the {\it effective action} ${\mathcal W}(g)$ is the
average over the configuration space of the exponentiated classical
action
$${\mathcal W}(g):= \langle {\mathcal A}(g)\rangle:= - \log \langle e^{-{\mathcal
    A}(g)}\rangle,$$ where $$\langle {\mathcal F} \rangle= \int_{{\rm Conf}(M)}
{\mathcal F}(\phi)\, {\mathcal D}\phi$$ 
is the average of ${\mathcal F}$
over the fields $\phi$ with respect
to some heuristic volume measure ${\mathcal D}\phi$ on
${\rm Conf}(M).$  The associated two-tensor $T_g({\mathcal W})$ is
interpreted as the {\it quantized stress-energy} momentum tensor and denoted
by
$$ T_g({\mathcal W})= T_g\left( \langle{\mathcal A}\rangle\right).$$ 
The local conformal anomaly 
associated to 
$d {\mathcal A}_g(f)= 
 2\left( {\rm tr}_g\left(T_g\left({\mathcal A}\right)\right),f\right)_g$
is defined  to be \cite{Du}
$${\rm tr}_gT_g\left( \langle{\mathcal A}\rangle\right)- \langle{\rm tr}_g
T_g({\mathcal A})\rangle.$$ 
If
the classical action is conformally invariant, as in string theory,
${\rm tr}_g
T_g({\mathcal A})=0$ and the local conformal anomaly reduces to 
${\rm tr}_gT_g\left( \langle{\mathcal
A}\rangle\right)$.

In general, the classical action is quadratic:
${\mathcal A}(g)(\phi)= \langle A_g\phi, \phi\rangle_g,$ 
where $\langle\cdot, \cdot\rangle_g$ is the inner product on the tensor fields
$\phi$ 
induced by the metric $g$.  $A_g$ is a
geometric differential operator, i.e.~an operator 
depending smoothly on the metric $g$ (via the
curvature, for example). For bosonic
strings, the fields are $\R^d$-valued smooth functions on a Riemann surface
$M$, and $A_g$ is the Laplace-Beltrami operator.
As pointed out in the introduction, even if 
${\mathcal A}(g)$ is conformally invariant, 
$A_g$ is usually
only conformally covariant.
\end{exm}

\subsection{Conformally covariant operators}

Given a vector bundle $E$ over a closed manifold $M$, we consider maps
$$
\Met(M)\to  Cl(M, E),\ \ 
g \mapsto  A_g.
$$

\begin{defn}
The operator $A_g\in Cl(M, E)$ associated to a Riemannian metric $g$ is 
{\it conformally covariant} of {\rm bidegree} $(a,b)$ if 
the pointwise
scaling of the metric $\bar g= e^{2f} g$, for $f\in \Ci(M, \R)$ yields
\begin{equation}\label{eq:confop}
 A_{\bar g}= e^{-b f} A_g e^{af}= e^{(a-b)f}\,A_g^\prime, 
\ \  {\rm for}\ A_g^\prime:=e^{-af} \, A_g\, e^{af},
\end{equation} 
for constants $a, b \in\R$.
\end{defn}

We survey known conformally covariant differential and pseudodifferntial
operators; more details are in Change \cite{C}.

\noindent {\bf Operators of order $1$.} (Hitchin \cite{H}) 
For $M^n$ spin,
the Dirac operator $\D_g:= \gamma^i \cdot \nabla^g_i$ 
is a conformally covariant operator of bidegree $\left(\frac{n-1}{2},
\frac{n+1}{2}\right)$. \\

\noindent {\bf Operators of order $2$.}  
If ${\rm dim}(M)=2$, the Laplace-Beltrami
operator $\Delta_g$
is  conformally covariant of bidegree $(0, 2)$.
It is well known that in dimension two
   \begin{equation}\label{eq:confscalar}
   R_{\bar g}= e^{-2f} \left(R_g+2\Delta_g f\right),
   \end{equation}
and by the Gauss-Bonnet theorem
   \begin{equation}\label{eq:GaussBonnet}
   \int_M R_g\ \dvol_g= 2\pi \chi(M),
   \end{equation}
with  the Euler characteristic $\chi(M)$ (much more than)
a conformal invariant. 

On a Riemannian manifold of dimension $n$,  
the Yamabe operator, also called the conformal Laplacian,
$$L_g:=  \Delta_g+c_n\, R_g,$$
is a conformally covariant operator  of bidegree $\left( \frac{n-2}{2},
 \frac{n+2}{2}\right)$,
where $R_g$ is the scalar curvature and $c_n:=\frac{n-2}{4(n-1)}$.


\noindent {\bf Operators of order $4$.} (Paneitz \cite{Pa,BO})
In  dimension $n$, the Paneitz operators 
$$P_g^n:= \tilde P_{g}^n+ (n-4) Q_{g}^n$$
are
conformally covariant scalar operators of bidegree $\left(\frac{n-4}{2},
\frac{n+4}{2}\right)$. 
Here $\tilde P_{g}^n:= \Delta_g^2$ \hfill\break
$ + d^* 
 \left( (n-2)\, J_g \, g-4  A_g\cdot \right) \, d$
with 
$$J_g:= \frac{R_g}{2(n-1)},  
A_g= \frac{Ric_g- \frac{R_g }{n}\, g}{n-2} +\frac{J_g}{n} g,$$
and $A_g\cdot$ the homomorphism on $T^*M$ given by
$\phi= (\phi_i)\mapsto
\left(A_g\right)_i^j\phi_j$,
and $Q_{g}^n:= \frac{n\, J_g^2- 4\vert A_g\vert^2 +2\Delta_g \,
J_g}{4}$ is Branson's
$Q$-curvature \cite{B1}, a local scalar invariant that
is a polynomial in the coefficients of the metric tensor and its inverse, the
scalar curvature and the Christoffel symbols.
Note that $A_g = \frac{1}{n}J_g\ g$ precisely when $g$ is Einstein.
\\ 

 The $Q$-curvature generalizes the scalar curvature $R_g$ in the following
 sense. On a $4$-manifold, we have
$$Q_{\bar g}^4= e^{-4f} \left(Q_g^4+\frac{1}{2}P_g^4 f\right)$$
(cf.~(\ref{eq:confscalar})), and
$\int_M Q_{g}^4 d{\rm vol}_g$ is a conformal invariant 
(cf.~(\ref{eq:GaussBonnet})), as is $\int_M Q_{g}^n {\rm dvol}_g$ in even
 dimensions \cite{B2}.


\noindent {\bf Operators of order $2k$.} (Graham, Jenne, Mason and Sparling \cite{GJMS})
Fix $k\in \Z^+$ 
and assume either $n$ is odd or $k\leq n$.  There are conformally
covariant (self-adjoint) scalar
differential operators $P_{g, k}^n$ 
of bidegree $\left(\frac{n-2k}{2}, \frac{n+2k}{2}\right)$ such that
the leading part of $P_{g, k}^n$ is $\Delta_g^k$ and such that
$P_{g, k}^n=\Delta_g^k$ on $\R^n$ with the Euclidean metric. 

$P_{g, k}^n$
generalizes $P_g^n$, since $P_g^n= P_{g, 2}^n$, and satisfies
$$P_{g, k}^n= \tilde P_{g}^n +\frac{n-2k}{2} Q_{g}^n$$
where $\tilde P_{g}^n= d^* S_g^n d$ for 
a natural differential operator $S_g^n$ on $1$-forms. 

Note that
$P_{g, k}^n$ has bidegree $(a,b)$ with $b-a=2k$ independent of
the dimension
 and in particular has bidegree $(0, 2k)$ in dimension $2k$.

\noindent {\bf Pseudodifferential Operators.} (Branson and Gover \cite{BG},
Petersen \cite{Pe})
  Peterson has
constructed $\pdo$s, $P_{g, k}^n, k\in \C,$ of order
$2{\rm Re}(k)$ and bidgree $((n-2k)/2, (n+2k)/2)$
on manifolds of dimension $n\geq 3$ with the property that
$P_{\bar g, k}^n - e^{-bf}P_{g,k}^n e^{af}$ is a smoothing operator.
Thus any conformal covariant built from the
total symbol of $P_{g, k}^n$ is a conformal covariant of $P_{g, k}^n$ itself.
The family $P_{g, k}^n$ contains the previously discovered
conformally covariant $\pdo$s 
associated to conformal boundary value problems \cite{BG}.

\subsection{A hirearchy of functionals and their conformal anomaly }

Since the known conformal invariants for conformally covariant operators $A_g$
$$\zeta_{A_g}(0)= {\rm tr}^{A_g}(I),\quad \log {\rm
det}_\zeta (A_g)= {\rm tr}^{A_g}\left(\log A_g\right), \quad \eta_{A_g}= {\rm
tr}^{ A_g} \left(A_g \, \vert A_g\vert^{-1}\right)$$ 
arise as weighted/canonical
 traces by (\ref{star1}), (\ref{star2}), (\ref{star3}),
it is natural to look for a
general prescription to derive conformal invariants from canonical traces.  

 Let $A_g\in Cl(M,E)$ 
be an operator associated to a Riemannian metric $g$ on $M$.  For $f\in
\Ci(M,\R)$, set $ g_t:= e^{2ft } g, t\in \R$, and set $A_t = A_{g_t}.$
We always assume that the map $g\mapsto A_g$ is smooth in the appropriate
topologies, so that $A_t$ is automatically a smooth curve in $\Cl(M,E).$ 

\begin{lem}
$A_g\in \Cl(M, E)$ 
is  conformally covariant of bidegree $(a, b)$ if and only if for all $f\in
\Ci(M, \R)$, 
\begin{equation}\label{eq:dbox}
\dot A_t=(a-b)\, f \,A_{t}-a\, [f, A_{ t}]. 
\end{equation}
\end{lem}

{\bf Proof:}
This follows from differentiating (\ref{eq:confop}) applied to the family
$g_t$.\\ ${}$ \hfill$\Box$ \\


\begin{thm}\label{thm:confanom} 
Let $A_g$ be a conformally covariant weight of bidegree $(a,b)$ and 
whose order
$\alpha$ and  spectral cut $\theta$ are independent of the metric
$g$. The meromorphic
map
$${\mathcal F}_h(g): z\mapsto  {\rm TR}\left( h(A_g) \, A_g^{-z}\right)$$ 
has conformal anomaly 
\begin{eqnarray}\label{eq:TRhbox}
 \delta_f {\rm TR}\left( h(A_g) \, A_g^{-z}\right)
&=& (a-b)\,
 {\rm TR}\left( f\,  h^\prime(A_g) \,
A_g^{-z+1}\right)\nonumber\\
&&\qquad -z(a-b) {\rm TR}\left(f\,  h(A_g) \, A_g^{-z}\right)
\end{eqnarray}
 as an identity of 
 meromorphic functions.

This is equivalent to the following system of equations. 
\begin{enumerate} 
\item[1.] The conformal anomaly of ${\rm res}(h(A_g))$:
\begin{equation}\label{eq:TRhresbox}
\delta_f {\rm res} \left( h(A_g) \right)
= (a-b)\, {\rm res}\left( f\, h^\prime(A_g) \, A_g\right).\end{equation}
\item[2.] The conformal anomaly of ${\rm tr}^{A_g}(h(A_g))$:
\begin{equation}\label{eq:TRh0box}
\delta_f {\rm tr}^{A_g} \left( h(A_g) \right)
= (a-b)\, {\rm tr}^{A_g}\left( f\, h^\prime(A_g) \, A_g\right)
+
\frac{b-a}{\alpha}\, {\rm res}(f\, h(A_g)), 
\end{equation}
\item[3.] The conformal anomaly of ${\rm tr}^{A_g}(h(A_g)\,\log^jA_g)$
 for  $j\in \Z^+$:
\begin{eqnarray}\label{eq:TRhjbox}
\lefteqn{\delta_f {\rm tr}^{A_g}\left( h(A_g)\, \log^jA_g\right)\ \ \ \ \ \ \
  \ \ \ \ \ \ \ \ }\nonumber\\
&=&(a-b)\,{\rm tr}^{A_g}\left(f\, h^\prime(A_g) \,A_g\, 
\log^jA_g\right)\\
&&\qquad +j\, (a-b)\, {\rm tr}^{A_g}\left( f\, h(A_g) \,
\log^{j-1}A_g\right).\nonumber
\end{eqnarray}
\end{enumerate}
\end{thm}

{\bf Proof:} Equations (\ref{eq:TRhbox}),
(\ref{eq:TRhresbox}),(\ref{eq:TRh0box}), (\ref{eq:TRhjbox}) follow from
(\ref{eq:TRh}), (\ref{eq:TRhresterm}),
(\ref{eq:TRh0term}), (\ref{eq:TRhjterm}), respectively.  In the computation,
we  use the cyclicity of TR on 
noninteger order operators, which eliminates the second term on the right hand
side of (\ref{eq:dbox}).  \hfill $\Box$\\

We collect these formulas for special choices of $h$.  We assume 
$A_g$ and hence $A_{\bar
  g}$ is invertible.  This allows us to ignore terms depending on 
the the kernel of $A_g$, which can be easily treated as in the proof of part 1
 below. 
All invariants and
covariants are understood to be conformal.

\begin{cor} \label{cor:confinv}
We have the following conformal anomalies for conformally covariant weights
$A_g$ of order $\alpha$:
\begin{enumerate}
\item[1.] Anomalies associated to $h\equiv 1$:
\begin{eqnarray}\label{h=1} 
\delta_f \zeta_{A_g}(0) &=&  \delta_f\tr^{A_g}(I) = 0,\nonumber\\
\delta_f \zeta^\prime_{A_g}(0) &=&
 - \delta_f\tr^{A_g}\left(\log\ A_g\right) = (a-b)\tr^{A_g}(f).
\end{eqnarray} 
Hence $\zeta_{A_g}(0) = -\frac{1}{\alpha}\res(\log A_g)$
is an invariant.
 $\zeta^\prime_{A_g}(0)$ has local anomaly
$-\frac{a-b}{\alpha}\res_x(\log A_g)$ and is an invariant whenever
$\res_x(\log A_g)$ vanishes for all $x\in M$.
\item[2.] Anomalies associated to $h(\lambda) = \lambda$:
\begin{eqnarray}
\delta_f\res(A_g) &=& (a-b)\res({f}A_g),\nonumber\\
 \delta_f\tr^{A_g}(A_g) &=& (a-b)\tr^{A_g}({f}A_g) +
\frac{b-a}{\alpha}\res({f}A_g), \label{h=lambda} \\
\delta_f\tr^{A_g}\left(A_g\log\ A_g\right) &=& 
(a-b)\tr^{A_g}({f}A_g\log A_g)  \nonumber\\
&&\qquad +(a-b)\tr^{A_g}({f}A_g).\label{h=lambda_2}
\end{eqnarray} 
$\res_x(A_g)$ is a pointwise covariant of weight $a-b$.
If $A_g$ is a differential operator, then $\res_x(A_g)$ vanishes 
and ${\rm tr}^{A_g}(A_g)=-\frac{1}{\alpha}\res(A_g\log A_g)$ has
local anomaly given by $\frac{b-a}{\alpha}\res_x(A_g\log A_g)$.
\item[3.] Anomalies associated to $h(\lambda) = \lambda^c, c\in \R$:  assuming
  $A_g$ is admissible for fixed $c$, we have
\begin{eqnarray}\label{h=lambdak} 
\delta_f\res(A_g^c) &=& c(a-b)\res(fA_g^c),\nonumber\\
\delta_f\tr^{A_g}(A_g^c) &=& c(a-b)\tr^{A_g}({f}A_g^c) 
+ \frac{b-a}{\alpha}\res({f}A_g^c).
\end{eqnarray} 
If $A_g^c$ is a differential operator (e.g.~if $A_g$ is differential and 
$c\in\Z^+$), then  $\res_x(A_g^c)$ vanishes 
and $\tr^{A_g}(A_g^c)= -\frac{1}{\alpha}\res(A_g^c \log A_g)$ 
has local anomaly given by
$\frac{c(b-a)}{\alpha}\res_x(A_g^c \log A_g)$.
\item[4.] Anomalies associated to $h(\lambda) = \lambda/|\lambda|$ and invertible
  $A_g$:  
\begin{equation}\label{h=sgn}
  \delta_f\tr^{A_g}(A_g/|A_g|)  = 
\frac{b-a}{\alpha}\res\left({f}\frac{A_g}{|A_g|}\right).
\end{equation}
\end{enumerate}
\end{cor}

{\bf Proof:} 
Much of the Corollary follows immediately from the Theorem.  In 1, 
the invariance of $\res_x(\log A_g)$ follows from (\ref{quick}). The last
statement follows from
 $\tr^{A_g}(f )=  -\frac{1}{\alpha} \res(f \log A_g)$ (\ref{eq:trQres}),
since multiplication by $f$ is a differential operator, and (\ref{star3}).
If Ker$(A)$ is nontrivial, $\zeta_{A_g}(0)$ is still a conformal invariant: by
 (\ref{star2}), (\ref{h=1}), the new terms
cause no trouble, as $\tr(\Pi_A) = {\rm dim}({\rm Ker}\ A)$ 
is a conformal invariant and 
$\res(I) = 0.$  

For the statement about $\res_x(A_g)$ in 2, 
if $\phi$ is a smooth function on $M$, then $\phi\cdot A_g$ is conformally
covariant if $A_g$ is. By (\ref{h=lambda}),
\begin{eqnarray*} \lefteqn{\delta_f \int_M \phi(x)
\left(\int_{S_x^*M} \tr_x \sigma_{-n}(A_g)(x,
\xi)d\xi\right) dx\ \ \ \ \ \ \ \ \ \ \ \ }\\
 &=&
\delta_f\res(\phi\cdot A_g) = (a-b)\res(f\cdot\phi\cdot A_g)\\
&=& (a-b) \int_M f(x)\phi(x)\left(\int_{S_x^*M} \tr_x \sigma_{-n}(A_g)(x,
\xi)d\xi\right) dx.
\end{eqnarray*}
Letting $\phi$ approach a delta function at $x$ and using the compactness of
$M$ to push this limit past
$\delta_f$ gives
$$\delta_f \res_x(A_g) = (a-b)f(x)\res_x(A_g).$$

The last statement in 2 follows as in the proof of 1 from
$\tr^{A_g}(A_g)=
-\frac{1}{\alpha} \res(A_g\log A_g)$.
In 3, the last statement follows from $\tr^{A_g}(f A_g^c)=
-\frac{1}{\alpha} \res(f A_g^c \log A_g).$
\hfill $\Box$\\

\begin{rem} The conformal anomaly in string theory boils down  to a finite
linear combination of local conformal anomalies $\frac{b-a}{\alpha}
\res_x(\log A_g)$ where the $A_g$ are Laplacians on forms. Their bidegree
involving the dimension of spacetime, so this local conformal anomaly
vanishes for a certain well chosen dimension.
\end{rem}

As stated in the introduction, the Corollary and the Laurent expansion of
Theorem \ref{thm:confanom} provide a natural
hierarchy among these invariants/covariants.
The most divergent term in the
 Laurent expansion is a conformal invariant; 
 if this global invariant vanishes in a particular case, 
 then the new ``most divergent''  term,
if it is of the form $\int_M {\mathcal I}(g, x) \dvol_g(x)$,
tends to give rise to a local conformal anomaly proportional to ${\mathcal I
}(g,x)$.
This is confirmed 
by the more refined  analysis for weights with nonnegative
leading order symbol, i.e.~weights with smoothing heat kernels.

\begin{lem}
 Let $A_g$ be a weight of order $\alpha$ with nonnegative leading symbol 
and let the heat kernel
for $A=A_g$ have the asymptotic expansion \cite{GS}
\begin{equation}\label{abc}
\tr_x e_{A}(\e,x,x)\sim \sum_{j=0}^\infty a_j(A,x) 
\e^{\frac {j-n}{\alpha}} + 
 \sum_{k=0}^\infty b_k(A,x)\e^k\log\e + \sum_{\ell=1}^\infty
 c_\ell(A,x)\e^\ell.\end{equation}
Then
\begin{equation}\label{coefs} \res_x(A^{k}) = \left\{\begin{array}{ll}
(-1)^{k+1}k!\alpha b_{k}(A,x),&  k\in \Z^{\geq 0},\\
  & \\
\frac{\alpha}{(-k-1)!}a_{n+\alpha k}(A, x),& k\in \Z^-,\end{array}\right.
\end{equation}
with the understanding that
 $a_{n+\alpha k} =0$ if $\alpha
k\not\in \Z.$   In particular, $\res_x(A) = \alpha b_{1}(A,x).$  

\end{lem}

The last sum in (\ref{abc})
appears only if $(j-n)/\alpha$ is never integral. In
 particular, this sum does not appear if $A_g$ is a differential operator.\\

{\bf Proof:}
Setting
$\zeta_A(z, x):= \omega_{KV}(A^{-z})(x) $
we have
$${\rm Res}_{z=-k} \zeta_A(z, x)= {\rm Res}_{z=0} \omega_{KV}(A^{-z+k})(x)=
-\frac{1}{\alpha'(0)}\res_x(A^k) = 
\frac{\res_x(A^k)}{\alpha},$$
(see (\ref{cutoffpatch}), (\ref{starstar})),
where we use the same symbol for an operator and its kernel and assume $A$ is
invertible for simplicity.
Let us compute this complex residue.
We have
\begin{eqnarray*}
\zeta_A(z,x) &=& \frac{1}{\Gamma(z)}\int_0^\infty t^{z-1} e_A(t,x,x) dt\\
&\sim & \frac{1}{\Gamma(z)}\int_0^1 t^{z-1}\left( \sum_{j=0}^\infty 
a_j(A,x)
t^{\frac {j-n}{\alpha}} +
 \sum_{k=0}^\infty  b_k(A,x)t^k\log t\right. \\
&&\qquad \left. +  \sum_{\ell=0}^\infty  c_\ell(A,x)t^\ell \right)dt
+ \frac{1}{\Gamma(z)}\int_1^\infty  t^{z-1}  e_A(t,x,x) dt.
\end{eqnarray*}
Since the last term is analytic in $z$, an easy integration on the first term
yields the result,
In particular, $\res_x(A) = \alpha b_{1}.$  
 \hfill $\Box$\\

To state the final theorem, 
let the kernel $ \tilde
e(\e,x,y)$ 
of
$Ae^{-\e |A|}$ have the asymptotic expansion
$$
\tr_x \tilde e(\e,x,x) \sim \sum_{j=0}^\infty \tilde
a_j(A,x)\e^{\frac{j-\alpha -n}{n}} + 
 \sum_{k=0}^\infty \tilde b_k(A,x)\e^k\log\e + \sum_{\ell=1}^\infty
 \tilde c_\ell(A,x)\e^\ell.
$$
Set
$a_j(A) = \int_M  a_j(A,x) \dvol_g(x)$, etc.

\begin{thm}\label{lastthm}
Let $A= A_g$ be a conformally covariant weight 
of bidegree $(a,b)$, with nonnegative leading order symbol, and 
whose order
$\alpha$ is independent of the metric
$g$.  
\begin{enumerate}
\item[1.] $a_n(A)+c_0(A)$  is a conformal invariant.  
\item[2.]  We have
\begin{eqnarray*}\delta_f\log\det{}_\zeta A &=&
-\delta_f(x) \zeta'_{A}(0)\\
&=& (b-a)\int_M f(x)\left( a_n(A,x)
+ c_0(A,x)\right) \dvol_g(x),
\end{eqnarray*}
so  $\log \det_\zeta(A)$ has local conformal
anomaly given by $(b-a)( a_n(A, x)+ c_0(A, x))$.
  In
particular, $\det_\zeta A$ is a conformal invariant if $A$ is a differential
operator and ${\rm dim}(M)$ is odd.
\item[3.] $b_1(A,x)$ is a 
a pointwise conformal covariant of weight $a-b.$
\item[4.]
$a_{n-\alpha}(A,x)$
is a pointwise conformal covariant of weight $b-a.$
\item[5.] We have
\begin{eqnarray*}
\lefteqn{\delta_f \tr^A(A)}\\
&=& - \delta_f(c_1(A)+a_{n+\alpha}(A))\\
&=& (b-a) \int_M
f(x) \left[ c_1(A, x)+ a_{n+\alpha}(A, x)- b_1(A, x)\right]\dvol_g(x),
\end{eqnarray*}
so that $\tr^A(A)$ has a local conformal anomaly given by
$$(b-a)\left(c_1(A, x)+ a_{n+\alpha}(A, x)- b_1(A, x)\right)$$
(with the understanding that
     $a_{n+\alpha}(A,x)  = 0$ if $\alpha\not\in\Z$, and $c_1(A,x) = 0$ if $A$
     is differential).
In particular, if $A$ is a differential operator, then $\tr^A (A)$
has a local  conformal anomaly given by $(b-a) a_{n+\alpha}(A, x)$,
and if $A$ is a noninteger order $\pdo$, it has a local conformal
anomaly given by $(b-a) c_1(A, x)$.
\item[6.]  
The results of 5 generalize to $\tr^A(A^k)$ for  $k\in \Z^+$, replacing
$c_1$ by $c_k$, $a_{n+\alpha}$ by $a_{n+\alpha k}$, $b_1$ by $b_k$, and
$(a-b)$ by $k(a-b)$.
%
\item [7.]$\delta_f\eta_{A}(0) =
  \frac{b-a}{\alpha}\res\left(f\frac{A}{|A|}\right)=  
-(b-a)\int_M f(x)\tilde a_n(A,x)$.  In particular,
  $\eta_A(0)$ is a conformal invariant if $n$ and $\alpha$ have opposite
  parity.  
\end{enumerate}
\end{thm}

{\bf Proof:} 1. This follows from the first point in
Corollary \ref{cor:confinv}
and the fact that $\zeta_{A_g}(0) =  a_n(A)+c_0(A)$.

2. This follows from (\ref{star3}), (\ref{h=1}), and the fact that
   $\tr^{A}({f} )
   = \int_M f(x) a_n(A,x).$  It is well known that the only the $\e^{k
   -\frac{{\rm dim}(M)}{2}}$ terms
in the   heat kernel asymptotics are nonzero for differential operators, so
   $a_n(A,x)=0$ in odd dimensions.

3.  This was shown in the second point of Corollary \ref{cor:confinv}.

4.  If $A$ is conformally covariant, so is $A^{-1}$.  The result now follows
    from 3 and the Lemma.

   5.   $\tr^{A}(A)
= \fp_{\e= 0} \tr(A e^{ -\e A})$ 
 is the finite part of $\tr(Ae^{-\e A}) = -\partial_{\e}
        \tr(e^{-\e A})$ as $\e\to 0,$  so $\tr^A(A) = -\int_M
        a_{n+\alpha}(A,x) + c_1(A,x).$  The
        first statement now follows from the second equation in
        (\ref{h=lambda}) and the fact that 
 $\frac{b-a}{\alpha}
\res(fA) = - b_1(A, x)$ (\ref{coefs}).
   If $A$ is a differential operator, then $c_1$ is
        replaced by $a_{n+\alpha}$ and $\res_x(fA) = 0$.  If $A$ is a
        non-integral order $\pdo$, then again $\res_x(fA) = 0.$  
   
   6.  By the first equation in 
   (\ref{h=lambdak}), $\delta_f\res(A^k) = k(a-b)\res(fA^k), k\in\Z^+.$
   The results for $b_k$ follow as in 5, using (\ref{coefs}).  
   We can use
   the second equation in (\ref{h=lambdak}) and
$\tr^A({f}A^k)  
= (-1)^k\fp_{\e=0}\partial_\e^k \tr( {f}e^{-\e A})$
   to prove the result for the $a$ and $c$ coefficients.

7.  The first equality 
follows from  (\ref{h=sgn}) and Remark 2.  
For the second equality, we have
$$\res(A/|A|) = \res_{s=0}\tr(A|A|^{-1}|A|^{-s}) = \res_{s=1}\tr(A|A|^{-s}).$$
Using the pointwise version of the
Mellin transform $A|A|^{-s} = \Gamma(s)^{-1}\int_0^\infty
t^{s-1}Ae^{ -t|A|} dt$ as in the Lemma, 
we get $\res\left(f\frac{A}{|A|}\right)=  
-(b-a)\int_M f(x)\tilde a_n(A,x)$.
The last
statement follows from a careful computation \cite[Prop.~3]{R}
 of the the residue of $A/|A|$.
Note that since $|A|$ has nonnegative symbol, this
restriction on the symbol of $A$ can be dropped here.
\hfill $\Box$

\begin{rem}  
(i) 1, 2, and 4 are known for the conformal Laplacian \cite{BO,PR}.
3 is new to
  our knowledge.  Related
  results for contact geometry are in Ponge \cite{Pon}.

(ii) The results above
  involving Wodzicki residues can be proved directly, where the cyclicity is
  valid for all order operators.
\end{rem}

\bibliographystyle{plain}

\end{document}